\documentclass[12pt,article]{amsart}

\usepackage{fullpage,amsfonts,amsmath,amssymb,amscd}
\input xypic

\theoremstyle{definition}
\newtheorem{defn}{Definition}[section]

\theoremstyle{plain}
\newtheorem{lem}[defn]{Lemma}
\newtheorem{prop}[defn]{Proposition}
\newtheorem{cor}[defn]{Corollary}
\newtheorem{thrm}[defn]{Theorem}
\newtheorem*{mainthrm}{Theorem}

\theoremstyle{remark}
\newtheorem{rem}[defn]{Remark}
\newtheorem{ex}[defn]{Examples}

\allowdisplaybreaks[3]

\newcommand{\N}{\mathbb{N}}
\newcommand{\Z}{\mathbb{Z}}
\newcommand{\C}{\mathbb{C}}

\title{The PI property of graded Hecke algebras}
\author{Katrin E Gehles}

\address{Alsenstrasse 10, 50679 Koeln, Germany}

\email{kgehles@web.de}

\begin{document}
\maketitle

\begin{abstract}
We show that graded Hecke algebras are PI algebras if and only if
they are finitely generated over their centres if and only if the
deformation parameters $t_{i}$ are zero for all $i=1,\ldots,N$. This
generalises a result for symplectic reflection algebras by Etingof -
Ginzburg and Brown - Gordon.
\end{abstract}

\section{Introduction}

Graded Hecke algebras were first defined by Drinfeld in \cite{Dri86}
and then studied in detail by Ram and Shepler in \cite{RS03}. Ram
and Shepler show that graded Hecke algebras are generalisations of
the graded affine Hecke algebras defined by Lusztig in \cite{Lus89}
for real reflection groups. Lusztig's work on graded affine Hecke
algebras was motivated by questions in group representation theory.
Drinfeld's construction of graded Hecke algebras is more general and
makes it possible to attach a graded Hecke algebra to any finite
subgroup of
$GL(V)$, not only the real reflection groups. 
In \cite{RS03} a full classification of the graded Hecke algebras
corresponding to complex reflection groups is achieved. Suprisingly
and disappointingly, there are complex reflection groups for which
no nontrivial graded Hecke algebras exist. This has inspired other
authors to look for further generalisations of graded Hecke
algebras, see \cite{SW06}.

Our work on graded Hecke algebras is inspired by their connection to
geometric questions. Graded Hecke algebras are deformations of skew
group algebras $S(V)\ast G$, where $V$ is a finite dimensional
vector space over $\C$ and $G$ a finite subgroup of $GL(V)$. The
centre of $S(V)\ast G$ is $S(V)^{G}$, the coordinate ring of the
orbit variety $V^{\ast}/G$. By studying $S(V)\ast G$ one hopes to
understand the $G$-equivariant geometry of $V^{\ast}$. If $V$ is a
symplectic vector space and the group $G$ preserves the symplectic
form, then the symplectic reflection algebras defined in \cite{EG02}
appear naturally as special cases of graded Hecke algebras. In this
special setting Etingof and Ginzburg are able to find smooth
deformations for some of the singular varieties $V^{\ast}/G$.

The purpose of this paper is to generalise the first steps in
\cite{EG02} to the more general setup of graded Hecke algebras. Our
results confirm a claim made by Etingof and Ginzburg in \cite[Remark
(ii), p. 246]{EG02}. The structure of symplectic reflection algebras
displays a dichotomy depending on the deformation parameter $t$.
Namely, a symplectic reflection algebra is finitely generated over
its centre if and only if it is a PI algebra if and only if $t=0$.
This is a result of \cite[Theorem 3.1]{EG02} and \cite[Proposition
7.2]{BG03}. Thus the obvious question to ask is whether graded Hecke
algebras display the same dichotomy in their behaviour depending on
specialisations of the deformation parameters. As it turns out there
is more than one deformation parameter which controls whether or not
graded Hecke algebras are PI. We denote these parameters by
$t_{1},\ldots,
t_{N}$. 
In the symplectic situation one can reduce to the case where $V$
is a symplectic vector space such that $V$ does not admit any
non-degenerate $G$-invariant subspaces. For such a vector space
the space of $G$-invariant skew-symmetric bilinear forms on $V$,
$((\bigwedge^{2}V)^{\ast})^{G}$, is one-dimensional, which gives
rise to the one parameter $t$. In general, however, the dimension
of the space $((\bigwedge^{2}V)^{\ast})^{G}$ can be greater than
one, say $N$, which leads to the appearance of $N$ deformation
parameters $t_{i}$. We show that
\begin{mainthrm} A graded Hecke algebra is
finitely generated over its centre if and only if it is a PI
algebra if and only if $t_{i}=0$ for all $i=1,\ldots, N$.
\end{mainthrm}

In Section \ref{grHA-definition} we begin by defining graded Hecke
algebras as in \cite{RS03} and derive some basic ring-theoretic
properties of these algebras. The definition of graded Hecke
algebras can also be motivated by
deformation theory, 
which is the content of Section \ref{grHA-deformation}. This
approach gives an explanation for the choice of the construction of
graded Hecke algebras and it will be crucial to proving our
subsequent results. Our work relies on the results in \cite{EG02}
and \cite[Proposition 7.2]{BG03} and the techniques developed in
\cite{EG02}. We modify their work to account for the fact that we
maintain a general setup and do not assume a symplectic structure.
Sections \ref{grHA-spherical} and \ref{grHA-preliminary} are devoted
to providing the details of these adjusted proofs. Finally, in
Section \ref{grHA-main} we prove our main theorem that tells us for
which values of the deformation parameters a graded
Hecke algebra has a big centre. 
As a corollary we deduce
the result mentioned above.

The work contained in this paper forms part of the author's PhD
thesis at the University of Glasgow. The author would like to thank
her supervisors K. A. Brown and I. Gordon for their advice and
support and C. Stroppel for many helpful comments. Throughout her
studies the author was supported by the Department of Mathematics of
the University of Glasgow and EPSRC.
\section{Definition and first properties}\label{grHA-definition}

Let $V$ be a finite dimensional vector space over $\C$ and $G$ a
finite subgroup of $GL(V)$. Denote by $\kappa:V\times V\rightarrow\C
G$ a skew-symmetric bilinear form. Let $T(V)=\C\oplus V\oplus (V
\otimes V)\oplus\cdots$ be the tensor algebra of $V$. The action of
$G$ on $V$ extends to an action of $G$ on $T(V)$ by $\C$-algebra
automorphisms. Construct the skew group algebra $T(V)\ast G$ and
define the following factor algebra
\[A:=\big(T(V)\ast G\big)/\langle[v,w]-\kappa(v,w)\,:\,v,w\in V\rangle,\]
where $[v,w]=vw-wv$.

The algebra $A$ is a positively filtered algebra, namely
\[F_{A}^{0}=\C G,\hspace{5mm} F_{A}^{1}=\C G +V\C G,\hspace{5mm}
F_{A}^{i}=(F_{A}^{1})^{i}\hspace{2mm}\text{for}\, i> 1.\] We form
the associated graded algebra $gr A$ of $A$ under this filtration.
From the relations in $A$ it is clear that there exists an
epimorphism $S(V)\ast G\twoheadrightarrow grA$. Here $S(V)$ denotes
the symmetric algebra of $V$, that is $S(V)=T(V)/\langle
[v,w]\,:\,v,w\in V\rangle$.

\begin{defn}\label{grHA:PBW}\cite[Section 1]{RS03} The algebra $A$ is called a \emph{graded Hecke algebra} if $S(V)\ast G\cong
grA$ as graded algebras. \end{defn}

Thus a graded Hecke algebra $A$ is isomorphic to $S(V)\otimes \C G$
as a graded vector space, which provides us with a PBW basis for
$A$. This PBW basis imposes necessary conditions on the form
$\kappa$. For example, $\kappa$ needs to be $G$-invariant in the
sense that $\kappa(g(v),g(w))=g\kappa(v,w)g^{-1}$, because otherwise
there exists a nontrivial linear relation in $A$ between elements of
$\C G$. In fact, there is a very precise description of the possible
choices for $\kappa$. To explain this we need to introduce more
notation.

Denote the centraliser of an element $g\in G$ by $Z_{G}(g)$. Recall
that a \emph{bireflection} in $G$ is an element $s\neq id\in G$ that
fixes a subspace of $V$ of codimension 2, that is
rank$_{V}(id-s)=2$. Write $V^{s}=\text{ker}(id-s)$ for the subspace
of $V$ fixed by $s$ and observe that $V\cong V^{s}\oplus V/V^{s}$,
where $V/V^{s}\cong\text{im}(id-s)$. Let $\mathcal{S}$ denote the
set of bireflections in $G$ and define
\begin{equation*}\label{grHA:subset-cond}
\mathcal{S}'=\{s\in \mathcal{S}\,|\, \forall g\in
Z_{G}(s),\hspace{0.2cm} \text{det}(g|_{V/V^{s}})=1\}
\end{equation*}
It follows directly from the definition that this set is closed
under conjugation. Let $S$ denote the normal subgroup of $G$
generated by $\mathcal{S}'$. Note that $S\subseteq SL(V)$ and,
therefore, $S$ does not contain reflections.

Let us construct specific skew-symmetric bilinear forms on $V$ as
follows. Fix $s\in\mathcal{S}'$. Since the space $\text{im}(id-s)$
is two-dimensional, there exists a unique - up to scalar
multiplication - nonzero skew-symmetric bilinear form on im$(id-s)$.
We can extend this form to all of $V$ by setting
$V^{s}=\text{ker}(id-s)$ to be its radical. Denote the form
constructed in this way by $\Omega_{s}$. Using this form as a
starting point we can define new forms: for $g\in G$ and $v,w\in V$,
\[\Omega_{g^{-1}sg}(v,w):=
\Omega_{s}(g(v),g(w)). \] With some easy calculations one can check
that the choice of the subset $\mathcal{S}'$ ensures that
$\Omega_{g^{-1}sg}$ is well-defined. Moreover, $\Omega_{g^{-1}sg}$
is indeed a nonzero skew-symmetric bilinear form on
im$(id-g^{-1}sg)$ with radical ker$(id-g^{-1}sg)$. Thus for a fixed
element $s\in \mathcal{S}'$ the forms corresponding to the elements
in the conjugacy class of $s$ are determined by $\Omega_{s}$.

Finally, let $\Omega$ be any skew-symmetric bilinear form on $V$
which is $G$-invariant, that is
$\Omega\in((\bigwedge^{2}V)^{\ast})^{G}$. The space
$((\bigwedge^{2}V)^{\ast})^{G}$ is a finite dimensional vector space
over $\C$. Let $\{b_{1},\ldots,b_{N}\}$ denote a basis for
$((\bigwedge^{2}V)^{\ast})^{G}$ over $\C$. Then
$\Omega=\sum_{i=1}^{N}t_{i}b_{i}$ for some $t_{i}\in \C$.\\

Now we are ready to state the crucial fact:

\begin{thrm}\label{grHA:def}\cite[Section 4]{Dri86}, \cite[Theorem 1.9]{RS03} With the above notation
the algebra $A$ is a graded Hecke algebra if and only if, for all
$v,w\in V$,
\[\kappa(v,w)=\Omega(v,w)\,id+\sum_{s\in
\mathcal{S}'}c_{s}\Omega_{s}(v,w)\, s,\] where
$\Omega(v,w)=\sum_{i=1}^{N}t_{i}b_{i}(v,w)$ for some $t_{i}\in\C$,
and the map $c:\mathcal{S}'\rightarrow\C$ given by $s\mapsto c_{s}$
is invariant under conjugation by $G$.
\end{thrm}

Therefore, a graded Hecke algebra is completely determined by the
choice of the complex values for $\{t_{i}\,|\,i=1\ldots,N\}$ and
$\{c_{s}\,|\,s\in \mathcal{S}'\}$. To express this fact in our
notation we henceforth denote a graded Hecke algebra by
$A_{\textbf{t},\textbf{c}}$, where $\textbf{t}$ denotes the
$N$-tuple of parameters $t_{i}$ and $\textbf{c}$ denotes the tuple
of parameters $\{c_{s}\,|\,s\in \mathcal{S}'\}$.

\begin{ex}
\
\begin{enumerate}
  \item If $t_{i}=0$ for all $i=1,\ldots,N$ and $c_{s}=0$ for all $s\in
\mathcal{S}'$, then $A_{\textbf{0},\textbf{0}}=S(V)\ast G$. This
follows directly from the defining relations of
$A_{\textbf{t},\textbf{c}}$. Thus the graded Hecke algebra
$A_{\textbf{t},\textbf{c}}$ is a deformation of the skew group
algebra $S(V)\ast G$.
  \item Let $V$ be a
symplectic vector space and $\omega$ the non-degenerate
skew-symmetric bilinear form on $V$. Suppose that the group $G$
preserves this form, that is $\omega(g(v),g(w))=\omega(v,w)$ for all
$g\in G$ and $v,w\in V$. Then symplectic reflection algebras, as
defined in \cite{EG02}, appear as special cases of graded Hecke
algebras. Namely, in the study of symplectic reflection algebras one
can reduce to the case where $V$ contains no non-degenerate
$G$-invariant subspaces. Under this assumption it can be shown that
$((\bigwedge^{2}V)^{\ast})^{G}=\C\omega$, see \cite[Section 2, p.
256]{EG02}. Thus without loss of generality we set $\Omega=\omega$.
Take a bireflection $s$. It is easy to see that the decomposition
$V=\text{ker}(id-s)\oplus\text{im}(id-s)$ is in
fact $\omega$-orthogonal. 
If we combine this with the fact that $\omega$ is non-degenerate on
$V$, we deduce that $\omega|_{\text{im}(id-s)}$ is non-degenerate.
Hence, without loss of generality we set
$\omega|_{\text{im}(id-s)}=\Omega_{s}$. Observe that in the
symplectic situation $\mathcal{S}=\mathcal{S}'$. But we have now
arrived at the definition of a symplectic reflection algebra, see
\cite[Theorem 1.3]{EG02}.
\end{enumerate}
\end{ex}



\begin{prop}\label{grHA:firstproperties}
Let $A_{\textbf{t},\textbf{c}}$ be a graded Hecke algebra.

(i) $A_{\textbf{t},\textbf{c}}= A_{\textbf{t},\textbf{c}}(S)\ast'
(G/S)$, where $A_{\textbf{t},\textbf{c}}(S)$ denotes the graded
Hecke algebra corresponding to $S$ instead of $G$, and $\ast'$
denotes some crossed product.

(ii) Let $\lambda\in\C^{\ast}$, then
$A_{\lambda\textbf{t},\lambda\textbf{c}}\cong
A_{\textbf{t},\textbf{c}}$.
\end{prop}

\begin{proof}
(i) 
Denote the ideal generated by
\[[v,w]=\kappa(v,w)=\Omega(v,w)\,id+\sum_{s\in \mathcal{S}'}c_{s}\Omega_{s}(v,w)s,\] for all $v,w\in V$, by
$I$. Then $A_{\textbf{t},\textbf{c}}=(T(V)\ast G)/I$ by definition.
By \cite[Lemma 1.3]{Pas89}, $(T(V)\ast G)/I=\big[(T(V)\ast S)\ast'
(G/S)\big]/I$.
Since $\kappa(v,w)\in \C S$ for all $v,w\in V$, the generators of
the ideal $I$ also generate an ideal of $T(V)\ast S$, which we
denote by $I_{S}$. We have $I=I_{S}\cdot(T(V)\ast G)=(T(V)\ast
G)\cdot I_{S}$ and $A_{\textbf{t},\textbf{c}}=\big[(T(V)\ast S)\ast'
(G/S)\big]/I=\big[(T(V)\ast S)/I_{S}\big]\ast'
(G/S)=A_{\textbf{t},\textbf{c}}(S)\ast' (G/S)$.

(ii) Define a map $A_{\textbf{t},\textbf{c}}\rightarrow
A_{\lambda\textbf{t}, \lambda\textbf{c}}$ by $x\mapsto
\sqrt{\lambda}x$ and $g\mapsto g$ for $x\in T(V), g\in G$.
\end{proof}

Because of the close connection between $A_{\textbf{t},\textbf{c}}$
and $A_{\textbf{t},\textbf{c}}(S)$ exhibited in this proposition we
will later restrict ourselves to the case when $G$ is generated by
the elements of $S$ that it contains, hence to the case $G=S$.

\begin{prop}\label{grHA:noeth,prime,gldim} Let $A_{\textbf{t},\textbf{c}}$ be a graded
Hecke algebra. Then $A_{\textbf{t},\textbf{c}}$ is noetherian, prime
and has finite global dimension.
\end{prop}

\begin{proof} We use the fact that $A_{\textbf{t},\textbf{c}}$ is filtered
such that $gr A_{\textbf{t},\textbf{c}}\cong S(V)\ast G$. By
\cite[Proposition 1.6.6, Theorem 1.6.9, Corollary 7.6.18]{MR87} all
the properties follow, if one can show that they hold for $S(V)\ast
G$. But $S(V)\ast G$ is noetherian and prime, because $G$ is a
finite group which acts faithfully on $S(V)$, see \cite[Proposition
1.6, Corollary 12.6]{Pas89}. Furthermore, since we are working over
$\C$, \cite[Theorem 7.5.6]{MR87} implies that the global dimension
of $S(V)\ast G$ is finite.
\end{proof}
\section{PBW deformation}\label{grHA-deformation}

In this section we show that the graded Hecke algebras
$A_{\textbf{t},\textbf{c}}$ are precisely a special kind of PBW
deformation of the skew group algebra $R:=S(V)\ast G$. This will
justify the shape of the form $\kappa$ which appears in Theorem
\ref{grHA:def}. Shepler and Witherspoon also consider graded Hecke
algebras as deformations of $R$ and determine the Hochschild
cocycles which appear, see \cite[Section 8]{SW06}.

The algebra $R=S(V)\ast G$ is naturally a positively graded algebra
with deg$\C G=0$. Let $R^{i}$ be the $i$th graded part of $R$, for
all $i\geq 0$, and note that each $R^{i}$ is a $\C G$-bimodule. Thus
we can view $R$ as a graded $\C G$-bimodule where the multiplication
in $R$ gives a $\C G$-bimodule map $R\otimes_{\C G} R\rightarrow R$.
From the relations in $A_{\textbf{t},\textbf{c}}$ it is clear that -
when constructing $A_{\textbf{t},\textbf{c}}$ as a deformation of
$R$ - we do not want to deform the relations in $\C G$. Thus to
ensure that we do not deform the degree zero part of $R$ we will use
the following definition in this section and only in this section.

\begin{defn} Let $B$ be a $\C G$-bimodule with a $\C G$-bimodule map $B\otimes_{\C G} B
\rightarrow B$. Then we call $B$ a \emph{$\C G$-algebra}.
\end{defn}

Therefore, subsequent maps will often be assumed to be $\C
G$-bimodule maps and we will frequently tensor over $\C G$ instead
of $\C$, which we will clarify by notation. In the literature
deformation theory usually happens over a field, but, as mentioned
in \cite[Section 2, p. 256]{EG02}, the theory explained in the
following and in particular the results of \cite{BG96} also
hold for $\C G$-algebras as defined above.\\

We recall the definition of a graded deformation of $R$. Suppose
$R_{h}$ is an associative unital algebra such that $R_{h}\cong
R\otimes_{\C}\C [h]$ as graded $\C$-vector spaces, where deg$h=1$.
Then $(R_{h}, \ast)$ is a graded deformation of $R$ if the
multiplication $\ast:R_{h}\times R_{h}\rightarrow R_{h}$ is a
$\C[h]$-linear map such that $r_{1}\ast r_{2}\equiv r_{1}r_{2}\,
\text{mod}\,hR_{h}$, for all $r_{1}, r_{2}\in R$. Thus
$R=R_{h}/hR_{h}$. In our situation we also require that $R_{h}$ is a
graded $\C G[h]$-bimodule and that $\ast$ is a $\C G[h]$-bimodule
map.

If $R_{h}$ is a graded deformation of $R$, then the multiplication
of two elements $r_{1},r_{2}\in R$ can be described by
\[r_{1}\ast r_{2}=r_{1}r_{2}+ \mu_{1}(r_{1},r_{2})\cdot h + \mu_{2}(r_{1},r_{2})\cdot h^{2} +
\cdots\, .\] The term $r_{1}r_{2}$ denotes the product in $R$ and
the maps $\mu_{i}:R\times R\rightarrow R$ are $\C G$-bimodule maps
of degree $-i$ with $i\in \N$. These maps completely determine the
multiplication in $R_{h}$ because of $\C G[h]$-linearity.

\begin{rem}\label{grHA:grdeform} Graded deformations have the following property: for all
$\lambda\in\C$ the factor algebra
$R_{h,\lambda}:=R_{h}/(h-\lambda)R_{h}$ is a filtered algebra such
that there is a canonical isomorphism $grR_{h,\lambda}\cong R$ as
algebras and also as $\C G$-algebras. The filtration on
$R_{h,\lambda}$ is induced by the filtration on $R_{h}$, which in
turn is derived from the grading on $R_{h}$.
\end{rem}

Let us now turn to the concept of a PBW deformation. We need to
introduce quadratic $\C G$-algebras first, see \cite[0.1 and
0.2]{BG96} for the definition of a quadratic algebra over a field.
In our case let $E$ denote a $\C G$-bimodule and let $T_{\C
G}(E)=\C\oplus E\oplus(E\otimes_{\C G} E)\oplus\ldots$ denote the
tensor $\C G$-algebra. Let $P$ be a subset of $\C G\oplus E \oplus
(E\otimes_{\C G} E)$ and also $\C G$-bimodule. If we denote the
ideal generated by $P$ by $I(P)$, then the algebra $Q(E,P):=T_{\C
G}(E)/I(P)$ is called a \emph{nonhomogeneous quadratic $\C
G$-algebra}. If $D$ is a subset of $E\otimes_{\C G} E$ and also a
$\C G$-bimodule, then we say that the algebra
$Q(E,D):=T_{\C G}(E)/I(D)$ is a \emph{quadratic $\C G$-algebra}. 

Suppose that $Q(E,P)$ is a nonhomogeneous quadratic $\C G$-algebra.
Then there exists a canonical quadratic $\C G$-algebra $Q(E,D)$
associated to $Q(E,P)$. Namely define $\pi:\C G\oplus E \oplus
(E\otimes_{\C G} E)\mapsto E\otimes_{\C G} E$ to be the projection
map and set $D=\pi(P)$. The $\C G$-algebra $T_{\C G}(E)$ is graded
with deg$\C G=0$ and deg$E=1$. This grading induces a filtration
$F^{\bullet}_{T_{\C G}(E)}$ on $T_{\C G}(E)$, which in turn induces
a filtration on $Q(E,P)$ via the surjection $p:T_{\C
G}(E)\twoheadrightarrow Q(E,P)$.
Namely, $F^{i}_{Q(E,P)}=p\,(F^{i}_{T_{\C G}(E)})=F^{i}_{T_{\C
G}(E)}/\big(F^{i}_{T_{\C G}(E)}\cap I(P)\big)$. The associated
graded $\C G$-algebra of $Q(E,P)$ under this filtration, denoted by
$grQ(E,P)$, is generated over $\C G$ by $p(E)$. Thus there exists a
surjective $\C G$-algebra map $T_{\C G}(E) \twoheadrightarrow
grQ(E,P)$. Since $\pi(P)=D$, we even have $\psi:
Q(E,D)\twoheadrightarrow grQ(E,P)$.

\begin{defn} \label{grHA:deform-PBW}\cite[Definition 0.3]{BG96} The nonhomogeneous
quadratic $\C G$-algebra $Q(E,P)$ is called a PBW deformation of
$Q(E,D)$ if $\psi$ is an isomorphism, that is if $Q(E,D)\cong
grQ(E,P)$.
\end{defn}

It is clear that any PBW deformation $Q(E,P)$ of $Q(E,D)$ must
satisfy the condition $P\cap F^{1}_{T_{\C G}(E)}=0$.
If this condition holds, the $\C G$-bimodule $P$ can be written
uniquely in terms of two $\C G$-bimodule maps $\alpha:D\mapsto E$
and $\beta:D\mapsto \C G$ as $P=\{d-\alpha(d)-\beta(d):d\in D\}$.\\


Let us see how the skew group algebra $R=S(V)\ast G$ and the graded
Hecke algebras $A_{\textbf{t},\textbf{c}}$ fit into this picture.
Take $E:=V\otimes_{\C}\C G$ and observe that $E$ is a free right $\C
G$-module by multiplication and a free left $\C G$-module by
$g(v\otimes g')=g(v)\otimes g\cdot g'$, for $v\in V$ and
$g,g'\in \C G$. 
Let $D\subseteq E\otimes_{\C G} E$ be the $\C$-span of
\[\{(v\otimes 1)\otimes_{\C G}(w\otimes g)-(w\otimes 1)\otimes_{\C G}(v\otimes
g)\},\] for $v,w\in V$ and $g\in \C G$. Note that $D$ is a $\C
G$-bimodule with the actions
\begin{flushleft}
$g'[(v\otimes 1)\otimes_{\C G}(w\otimes g)-(w\otimes 1)\otimes_{\C
G}(v\otimes g)]=$
\end{flushleft}
\begin{flushright}
$(g'(v)\otimes 1)\otimes_{\C G}(g'(w)\otimes g'g)-(g'(w)\otimes
1)\otimes_{\C G}(g'(v)\otimes g'g),$
\end{flushright}
\begin{center}
$[(v\otimes 1)\otimes_{\C G}(w\otimes g)-(w\otimes 1)\otimes_{\C
G}(v\otimes g)]g'=(v\otimes 1)\otimes_{\C G}(w\otimes gg')-(w\otimes
1)\otimes_{\C G}(v\otimes gg').$
\end{center}
We claim that $S(V)\ast G\cong Q(E,D)$ as $\C G$-algebras. To see
this we construct an isomorphism $\theta:T_{\C G}(E)\rightarrow
T(V)\ast G$, where $T(V)$ denotes the usual tensor algebra of $V$
over $\C$, as the direct sum of the maps $\theta_{p}: E^{\otimes_{\C
G} p}\rightarrow V^{\otimes p}\otimes \C G$, given by
\[\theta_{p}[(v_{1}\otimes g_{1})\otimes_{\C G}\cdots\otimes_{\C G} (v_{p}\otimes
g_{p})]=\big(v_{1}\otimes h_{1}(v_{2})\otimes
h_{2}(v_{3})\otimes\cdots\otimes h_{p-1}(v_{p})\big)\otimes h_{p},\]
where $v_{i}\in V$, $g_{i}\in \C G$ and $h_{i}:=g_{1}\cdots g_{i}$
for $i=1,\ldots, p$. 
Then $\theta_{2}(D)= C\otimes \C G$, where $C$ denotes the space of
commutators in $V^{\otimes 2}$. Thus $Q(E,D)=T_{\C G}(E)/I(D)\cong
S(V)\ast
G=R$ is a quadratic $\C G$-algebra.\\

\begin{lem} Let $A_{\textbf{t},\textbf{c}}$ be a graded Hecke
algebra and let $Q(E,D)\cong S(V)\ast G$ with $E$ and $D$ defined as
above. Then there exists $P\subseteq\C G\oplus E \oplus
(E\otimes_{\C G} E)$ such that $A_{\textbf{t},\textbf{c}}=Q(E,P)$, a
nonhomogeneous quadratic $\C G$-algebra. Moreover, the quadratic $\C
G$-algebra associated to $Q(E,P)$ under the projection $\pi:\C
G\oplus E \oplus (E\otimes_{\C G} E)\mapsto E\otimes_{\C G} E$ is
$S(V)\ast G$ and so $A_{\textbf{t},\textbf{c}}$ is a PBW deformation
of $S(V)\ast G$.
\end{lem}

\begin{proof} A graded Hecke algebra is defined as
$A_{\textbf{t},\textbf{c}}=(T(V)\ast G)/I$, where $I$ denotes the
ideal generated by $[v_{1},v_{2}]-\kappa(v_{1},v_{2})$, for all
$v_{1},v_{2}\in V$. Equivalently one could choose as generators of
the ideal $I$ the elements
$\big([v_{1},v_{2}]-\kappa(v_{1},v_{2})\big)g$, for all
$v_{1},v_{2}\in V$ and $g\in\C G$. We have seen on the previous page
that $T_{\C G}(E)\cong T(V)\ast G$ via an isomorphism which we
labelled $\theta$. The map $\theta$ induces an isomorphism
$A_{\textbf{t},\textbf{c}}\cong T_{\C G}(E)/I(P)=Q(E,P)$, where $P$
is the $\C$-span of
\[\{(v_{1}\otimes 1)\otimes_{\C G}(v_{2}\otimes g)-(v_{2}\otimes
1)\otimes_{\C G}(v_{1}\otimes g)-\kappa(v_{1},v_{2})g\},\] for
$v_{1},v_{2}\in V,g\in \C G$.
Note that we can extend the $\C G$-action on the subset $D$ to make
$P$ into a $\C G$-bimodule, because $\kappa$ is $G$-invariant. The
quadratic $\C G$-algebra naturally associated to
$A_{\textbf{t},\textbf{c}}\cong Q(E,P)$ is clearly $S(V)\ast G\cong
Q(E,D)$. Now, by definition $grA\cong S(V)\ast G$, and hence graded
Hecke algebras are PBW deformations of $S(V)\ast G$. Therefore, we
can write $P$ in terms of $\C G$-bimodule maps $\alpha$ and $\beta$.
Namely, let $\alpha=0$ and $\beta[(v_{1}\otimes 1)\otimes_{\C
G}(v_{2}\otimes g)-(v_{2}\otimes 1)\otimes_{\C G}(v_{1}\otimes
g)]=\kappa(v_{1},v_{2})g$.
\end{proof}


We now want to show that, conversely, all PBW deformations of
$S(V)\ast G$ with certain properties are graded Hecke algebras.
Observe that the $\C G$-algebra $S(V)\ast G$ is Koszul, which can be
seen from tensoring the Koszul resolution of the trivial
$S(V)$-module $\C$ on the right by $\C G$. We have the following
result


\begin{thrm}\label{grHA:deform-cond}\cite[Lemma 0.4, Lemma 3.3, Theorem 4.1]{BG96}
Let $Q(E,D)$ be a quadratic Koszul $\C G$-algebra, where $E$ is a
free $\C G$-module from either side. Assume that we are given
$Q(E,P)$ in terms of $\C G$-bimodule maps $\alpha:D\mapsto E$,
$\beta:D\mapsto \C G$, and $P=\{d-\alpha(d)-\beta(d):d\in D\}$. Then
$Q(E,D)\mapsto grQ(E,P)$ is an isomorphism if and only if the
following are satisfied

(i) $\alpha\otimes_{\C G} id-id\otimes_{\C G}\alpha$ has image in
$D$,

(ii) $\alpha\circ(\alpha\otimes_{\C G}id-id\otimes_{\C
G}\alpha)=id\otimes_{\C G}\beta-\beta\otimes_{\C G}id$,

(iii)
$\beta\circ(\alpha\otimes_{\C G}id-id\otimes_{\C G}\alpha)=0$,\\
where the domain of all these maps is $(D\otimes_{\C G}
E)\cap(E\otimes_{\C G} D)$.
\end{thrm}

\begin{rem}\label{grHA:proof-deform} In fact, in \cite[Theorem 4.1]{BG96} it is proved that the conditions (i) - (iii) on
the maps $\alpha$ and $\beta$ and the fact that $Q(E,D)$ is Koszul
allow one to construct a graded deformation
$\big(Q(E,D)_{h},\ast\big)$. Let
$Q(E,D)_{h,1}=Q(E,D)_{h}/\big((h-1)Q(E,D)_{h}\big)$. One then
obtains $\C G$-bimodule maps
\[Q(E,D)\xrightarrow{\psi} grQ(E,P)\xrightarrow{\rho} gr Q(E,D)_{h,1}\xrightarrow{\varphi} Q(E,D),\]
where $\psi$ is the natural surjection from above and, by Remark
\ref{grHA:grdeform} in this section, $\varphi$ is an isomorphism.
The map $\rho$ comes from the $\C G$-bimodule map that includes $E$
in $Q(E,D)_{h}$ and then projects onto $Q(E,D)_{h,1}$. This map
extends uniquely to an algebra and $\C G$-bimodule map $T_{\C
G}(E)\rightarrow Q(E,D)_{h,1}$, which factors through $Q(E,P)$. The
map $\rho$ is then the associated graded map of this map
$Q(E,P)\rightarrow Q(E,D)_{h,1}$. Finally, one checks that the
composition $\varphi\circ\rho\circ\psi$ is the identity map on
elements of degrees zero and one in $Q(E,D)$. Since $Q(E,D)$ is
generated by those elements, the composition is just the identity
map which implies that $Q(E,D)\cong grQ(E,P)$ and $grQ(E,P)\cong
grQ(E,D)_{h,1}$.
\end{rem}


\begin{cor}\label{grHA:deform}
Let $Q(E,D)\cong S(V)\ast G$ with $E$ and $D$ as before. Suppose
$Q(E,P)$ is a nonhomogeneous quadratic $\C G$-algebra and a PBW
deformation of $Q(E,D)$. Then the $\C G$-bimodule $P$ is given by
$\C G$-bimodule maps $\alpha:D\mapsto E$, $\beta:D\mapsto \C G$,
such that $P=\{d-\alpha(d)-\beta(d):d\in D\}$. Assume that
$\alpha=0$. Then $Q(E,P)$ is a graded Hecke algebra.
\end{cor}

\begin{proof}
See \cite[p.257]{EG02} for the first part of this proof. More
details are taken from \cite{Gor05} and are given here for the
reader's convenience.

Since $\beta$ is a right $\C G$-module map, it is determined by
some antisymmetric mapping, say $\kappa:V\times V\rightarrow \C G$.
Namely, $\beta[(v\otimes 1)\otimes_{\C G}(w\otimes 1)-(w\otimes
1)\otimes_{\C G}(v\otimes 1)]=\kappa(v,w)$.
The fact that $\beta$ is a $\C G$-bimodule map translates into
$\kappa$ being $G$-invariant in the sense that
$\kappa(g(v),g(w))=g\kappa(v,w)g^{-1}$ for all $g\in G$.
Only condition (ii) of the previous theorem is non-trivial and it
reduces to $0=id\otimes_{\C G}\beta-\beta\otimes_{\C G}id$ on
$(D\otimes_{\C G} E)\cap(E\otimes_{\C G} D)$. We can use the
isomorphism $\theta:T_{\C G}(E)\rightarrow T(V)\ast G$ again to
identify $D$ with $C\otimes \C G$, where $C$ denotes the space of
commutators in $V^{\otimes 2}$.
Since $\beta$ is a right $\C G$-module map, condition (ii) reduces
to $id\otimes\kappa-\kappa\otimes id=0$ on $(C\otimes V)\cap
(V\otimes C)$.
In \cite[p.257]{EG02} it is proved that this condition implies that,
for all $v,w\in V$, we have
$\kappa(v,w)=a_{1}(v,w)\,id+\sum_{s\in\mathcal{S}}a_{s}(v,w)s$, for
some skew-symmetric bilinear forms $a_{1},a_{s}:V\times
V\rightarrow\C$. Furthermore, it follows that $a_{1}\in ((\bigwedge
V)^{\ast})^{G}$ and that $V^{s}\subseteq\text{ker}a_{s}$ for all
$s\in\mathcal{S}$.

Now suppose $s\in \mathcal{S}$. Then the form $a_{s}$ must be
proportional to the skew-symmetric bilinear form $\Omega_{s}$ on
$V$, which we constructed in Section \ref{grHA-definition}. Thus
assume without loss of generality that $a_{s}=\Omega_{s}$ for all
$s\in\mathcal{S}$. Denote a basis of the 2-dimensional subspace
$\text{im}(id-s)$ by $v_{1},v_{2}\in V$ and extend this to a basis
$v_{1},\ldots,v_{n}$ of $V$. Suppose $g\in Z_{G}(s)$, then the
$G$-invariance of $\kappa$ implies that
$\Omega_{s}(v,w)=\Omega_{s}(g(v),g(w))$ for all $v,w\in V$. Since
$g(V^{s})=V^{s}$, this reduces to the statement
$\Omega_{s}(v_{1},v_{2})=\Omega_{s}(g(v_{1}),g(v_{2}))$. Write
$g(v_{1})=av_{1}+cv_{2}$ and $g(v_{2})=bv_{1}+dv_{2}$ for some
$a,b,c,d\in\C$. Then
\[\Omega_{s}(g(v_{1}),g(v_{2}))=\Omega(av_{1}+cv_{2},bv_{1}+dv_{2})=(ad-bc)\Omega_{s}(v_{1},v_{2}),\]
which implies $ad-bc=1$. Hence for all $g\in Z_{G}(s)$ we have
det$(g|_{V/V^{s}})=1$, so $s\in \mathcal{S}'$. Thus $\kappa$ must
have precisely the form described in Theorem \ref{grHA:def}.
\end{proof}

The fact that $A_{\textbf{t},\textbf{c}}$ is a PBW deformation of
$Q(E,D)\cong S(V)\ast G=R$ will become crucial in the proof of our
main theorem. In Remark \ref{grHA:proof-deform} in this section we
observed that there must exist a graded deformation $(R_{h},\ast)$
of $R$ such that
$R_{h}/\big((h-1)R_{h}\big)=A_{\textbf{t},\textbf{c}}$. For $v,w\in
V\subset R$ we have $v\ast w = vw +\mu_{1}(v,w)\cdot
h+\mu_{2}(v,w)\cdot h^{2}+\ldots$ in $R_{h}$, for some $\C
G$-bimodule maps $\mu_{i}:R\times R\rightarrow R$ of degree $-i$.
Since $v\ast w$ has degree $2$, we must have $\mu_{i}(v,w)=0$ for
$i>2$, hence $v\ast w = vw +\mu_{1}(v,w)\cdot h+\mu_{2}(v,w)\cdot
h^{2}$. Moreover, in the factor $R_{h}/\big((h-1)R_{h}\big)$ we have
$v\ast w \equiv vw +\mu_{1}(v,w)+\mu_{2}(v,w)$. On the other hand,
$R_{h}/\big((h-1)R_{h}\big)=A_{\textbf{t},\textbf{c}}$ and,
therefore, $v\ast w=w\ast v+\kappa(v,w)$. We deduce that
$[\mu_{1}(v,w)-\mu_{1}(w,v)]+[\mu_{2}(v,w)-\mu_{2}(w,v)]=\kappa(v,w)$.
Thus we must have $\mu_{1}(v,w)-\mu_{1}(w,v)=0$, because
deg$\mu_{1}(v,w)=\text{deg}\mu_{1}(w,v)=1$ and $\kappa(v,w)\in \C G$
has degree zero. In summary, for all $v,w\in V$, we have
\[\mu_{2}(v,w)-\mu_{2}(w,v)=\kappa(v,w)=\Big[\sum_{i=1}^{N}t_{i}b_{i}(v,w)\Big]\cdot id+\sum_{s\in
\mathcal{S}'}c_{s}\Omega_{s}(v,w)s.\] Similarly, one can see that
for all $p,p'\in S(V)$ the difference $\mu_{2}(p,p')-\mu_{2}(p',p)$
depends linearly on the parameters in $\textbf{t}$ and $\textbf{c}$.
\section{The spherical subalgebra}\label{grHA-spherical}

Let $A_{\textbf{t},\textbf{c}}$ denote a graded Hecke algebra as
defined in Section \ref{grHA-definition}. Recall that the
symmetrizing idempotent in $\C G\subseteq A_{\textbf{t},\textbf{c}}$
is given by $e=\frac{1}{|G|}\sum_{g\in G}g$. The \emph{spherical
subalgebra} of $A_{\textbf{t},\textbf{c}}$ is defined as
$eA_{\textbf{t},\textbf{c}}e$. It is easy to see that the filtration
on $A_{\textbf{t},\textbf{c}}$ intersects with the spherical
subalgebra to induce a filtration
$F^{\bullet}_{eA_{\textbf{t},\textbf{c}}e}$ on
$eA_{\textbf{t},\textbf{c}}e$. We have graded algebra isomorphisms
\[gr(eA_{\textbf{t},\textbf{c}}e)=e(gr A_{\textbf{t},\textbf{c}})e\cong
e(S(V)\ast G)e\cong S(V)^{G},\] where the inverse of the last
isomorphism is given by the map $p\mapsto pe$ for $p\in S(V)^{G}$.
Observe that $S(V)^{G}=Z\big(S(V)\ast
G\big)=Z(A_{\textbf{0},\textbf{0}})$.\\

The space $A_{\textbf{t},\textbf{c}}e$ has a left
$A_{\textbf{t},\textbf{c}}$-module structure and a right
$eA_{\textbf{t},\textbf{c}}e$-module structure, both given by
multiplication. Again the filtration of $A_{\textbf{t},\textbf{c}}$
induces a filtration $F^{\bullet}_{A_{\textbf{t},\textbf{c}}e}$ on
the module $A_{\textbf{t},\textbf{c}}e$. We have
$gr(A_{\textbf{t},\textbf{c}}e)\cong S(V)\cong
A_{\textbf{0},\textbf{0}}e$, which can be deduced by using the same
isomorphisms as above for $gr(eA_{\textbf{t},\textbf{c}}e)$.

\begin{lem}\label{grHA:spherical-domain}
\ \\
\indent(i) $eA_{\textbf{t},\textbf{c}}e$ is a finitely generated
$\C$-algebra and a noetherian domain.

(ii) $A_{\textbf{t},\textbf{c}}e$ is finitely generated as right
$eA_{\textbf{t},\textbf{c}}e$-module.
\end{lem}

\begin{proof}
This follows directly from associated graded techniques as in
\cite[Lemma 7.6.11]{MR87} and the Hilbert-Noether Theorem, see
\cite[Theorem 1.3.1]{Ben93}.
\end{proof}

Recall that $S$ denotes the subgroup of $G$ generated by the
elements of the set $\mathcal{S}'$, as defined in Section
\ref{grHA-definition}.

\begin{lem}\label{grHA:End}
Assume that $G=S$. Then
$A_{\textbf{t},\textbf{c}}\cong\emph{End}_{eA_{\textbf{t},\textbf{c}}e}(A_{\textbf{t},\textbf{c}}e)$
as algebras.
\end{lem}

\begin{proof} In large parts we use the proof of \cite[Theorem 1.5
(iv)]{EG02}. For the reader's convenience we give the full details.

Left multiplication by elements of $A_{\textbf{t},\textbf{c}}$
induces an algebra map $\eta: A_{\textbf{t},\textbf{c}}\rightarrow
\text{End}_{eA_{\textbf{t},\textbf{c}}e}(A_{\textbf{t},\textbf{c}}e)$
by $a\mapsto \big(l_{a}:a'e\mapsto aa'e\big)$, for $a,a'\in
A_{\textbf{t},\textbf{c}}$. This map is in fact filtration
preserving, where a filtration on
$\text{End}_{eA_{\textbf{t},\textbf{c}}e}(A_{\textbf{t},\textbf{c}}e)$
is defined as follows. Denote the generators of
$gr(A_{\textbf{t},\textbf{c}}e)$ as
$gr(eA_{\textbf{t},\textbf{c}}e)$-module by
$\overline{u}_{1},\ldots,\overline{u}_{n}$ and let
deg$(\overline{u}_{i})=d_{i}$. Then $A_{\textbf{t},\textbf{c}}e$ is
generated as $eA_{\textbf{t},\textbf{c}}e$-module by representatives
of the $\overline{u}_{1},\ldots,\overline{u}_{n}$ denoted by
$u_{i}\in A_{\textbf{t},\textbf{c}}e$, see the proof of \cite[Lemma
7.6.11]{MR87}. Now take an element $f\in
\text{End}_{eA_{\textbf{t},\textbf{c}}e}(A_{\textbf{t},\textbf{c}}e)$.
We can find $m\in\Z$ such that $f(u_{i})\in
F_{A_{\textbf{t},\textbf{c}}e}^{d_{i}+m}$ for all $i=1,\ldots,n$.
Therefore, $f(F^{j}_{A_{\textbf{t},\textbf{c}}e})\subseteq
F^{j+m}_{A_{\textbf{t},\textbf{c}}e}$ for all $j\geq 0$. Thus we
have an increasing $\Z$-filtration on
$\text{End}_{eA_{\textbf{t},\textbf{c}}e}(A_{\textbf{t},\textbf{c}}e)$:
\[
F^{m}_{\text{End}}=\{f\in\text{End}_{eA_{\textbf{t},\textbf{c}}e}(A_{\textbf{t},\textbf{c}}e)\,|\,
f(F^{n}_{A_{\textbf{t},\textbf{c}}e})\subseteq
F^{n+m}_{A_{\textbf{t},\textbf{c}}e}\hspace{3mm} \forall\,
n\in\Z\}.\] Since $\eta$ is filtration preserving, we can construct
the algebra map $gr(\eta)$. It now suffices to show that $gr(\eta)$
is an algebra isomorphism, see \cite[Corollary 7.6.14]{MR87}.
To this end we consider the composite
\[ grA_{\textbf{t},\textbf{c}}\xrightarrow{gr(\eta)}
gr\big(\text{End}_{eA_{\textbf{t},\textbf{c}}e}(A_{\textbf{t},\textbf{c}}e)\big)\xrightarrow{j}
\text{End}_{gr(eA_{\textbf{t},\textbf{c}}e)}\big(gr(A_{\textbf{t},\textbf{c}}e)\big),\]
where the map $j$ is given by $f+F^{m-1}_{\text{End}}\mapsto
\Big[ae+F^{k-1}_{A_{\textbf{t},\textbf{c}}e}\mapsto
f(ae)+F^{k+m-1}_{A_{\textbf{t},\textbf{c}}e}\Big]$, for $a\in
F^{k}$. The map $j$ is clearly injective
and we have reduced the problem to showing that the composite
$j\circ gr(\eta):S(V)\ast G\rightarrow
\text{End}_{S(V)^{G}}\big(S(V)\big)$ is an algebra isomorphism. Let
us first show that $j\circ gr(\eta)$ is injective. Observe that
$S(V)^{G}$ is central both in $S(V)\ast G$ and
$\text{End}_{S(V)^{G}}\big(S(V)\big)$. We tensor on the left with
the quotient field of $S(V)^{G}$, $Q\big(S(V)^{G}\big)$:
\[Q\big(S(V)^{G}\big)\otimes_{S(V)^{G}}\big(S(V)\ast
G\big)\xrightarrow{id\otimes (j\circ gr(\eta))}
Q\big(S(V)^{G}\big)\otimes_{S(V)^{G}}\text{End}_{S(V)^{G}}\big(S(V)\big).\]
Let us show that $\varphi:=id\otimes \big(j\circ gr(\eta)\big)$ is
an algebra isomorphism. We have the following isomorphisms as
$S(V)^{G}$-modules, see \cite[Lemma 2.4, Proposition 2.10]{Eis95},
which imply algebra isomorphism:
\[Q\big(S(V)^{G}\big)\otimes_{S(V)^{G}}\big(S(V)\ast
G\big)\hspace{2mm} \cong\hspace{2mm}
\Big[Q\big(S(V)^{G}\big)\otimes_{S(V)^{G}}S(V)\Big]\ast
G\hspace{2mm} \cong\hspace{2mm}
Q(S(V))\ast G,\] 
\begin{eqnarray*}
Q\big(S(V)^{G}\big)\otimes_{S(V)^{G}}\text{End}_{S(V)^{G}}(S(V))&
\cong &
\text{End}_{Q(S(V)^{G})}\Big[Q\big(S(V)^{G}\big)\otimes_{S(V)^{G}}S(V)\Big]\\
\ & \cong & \text{End}_{Q(S(V)^{G})}\big[Q(S(V))\big],\\
\end{eqnarray*}
where $Q\big(S(V)^{G}\big)\otimes_{S(V)^{G}}S(V)\cong Q(S(V))$,
because $S(V)$ is a finitely generated $S(V)^{G}$-module. 
The map $\varphi$ is given by $\sum_{g\in G}p_{g}g\mapsto
\big[x\mapsto \sum_{g\in G}p_{g}\cdot g(x)\big]$, for $p_{g}\in
Q(S(V))$. First observe that $\varphi$ is clearly not the zero map.
Then note that $Q(S(V))\ast G$ is a simple ring, since $Q(S(V))$ is
a simple ring and $G$ acts faithfully on $Q(S(V))$, see
\cite[Proposition 7.8.12]{MR87}. Thus ker$\varphi=0$. Now count the
dimensions of the $Q\big(S(V)^{G}\big)$-vector spaces on each side
of the map $\varphi$. We have dim$_{Q(S(V)^{G})}\big[Q(S(V))\ast
G\big]=|G|^{2}=\text{dim}_{Q(S(V)^{G})}\Big[\text{End}_{Q(S(V)^{G})}\big[Q(S(V))\big]\Big]$,
since $Q(S(V))$ is a Galois extension of $Q(S(V)^{G})$ and
$\big[Q(S(V)):Q(S(V)^{G})\big]=|G|$.

The fact that $\varphi=id\otimes \big(j\circ gr(\eta)\big)$ is an
isomorphism now implies that $j\circ gr(\eta)$ is injective, because
of the following commutative diagram:
\[\xymatrix
{Q(S(V))\ast G \ar[r]^>>>>{\varphi} & \text{End}_{Q(S(V)^{G})}\big[Q(S(V))\big]\\
S(V)\ast G \ar[u] \ar[r]^>>>>>>>>{j\circ gr(\eta)} &
\text{End}_{S(V)^{G}}\big(S(V)\big) \ar[u]\\}\] where the vertical
map on the left is an embedding, since the elements of $S(V)$ are
nonzero divisors of $S(V)\ast G$. 

It remains to show that $j\circ gr(\eta)$ is surjective, hence that
im$\big(j\circ gr(\eta)\big)=\text{End}_{S(V)^{G}}\big(S(V)\big)$.
To do this we
will use the concept of a maximal order. 
By \cite[Theorem 4.6]{Mar95}, $S(V)\ast G$ is a maximal order in its
quotient field if and only if $G$ does not contain reflections in
its action on $S(V)$. But we have assumed that $G=S\subseteq SL(V)$
and so $S$ does not contain reflections.
Moreover, the map $\varphi$ shows that the quotient rings of
$S(V)\ast G$ and $\text{End}_{S(V)^{G}}\big(S(V)\big)$ are
isomorphic. Note that the quotient ring of
$\text{End}_{S(V)^{G}}\big(S(V)\big)$ is indeed
$\text{End}_{Q(S(V)^{G})}\big[Q(S(V))\big]$. 

Now we use the commutative diagram from above again. As $S(V)\ast G$
is a maximal order in its quotient ring $Q(S(V))\ast G$ and $j\circ
gr(\eta)$ is injective, im$\big(j\circ gr(\eta)\big)$ is also a
maximal order in the quotient ring
$\text{End}_{Q(S(V)^{G})}\big[Q(S(V))\big]$. But
$\text{End}_{S(V)^{G}}\big(S(V)\big)\supseteq S(V)\ast G$ via the
embedding $j\circ gr(\eta)$ and
$\text{End}_{S(V)^{G}}\big(S(V)\big)$ is finitely generated over
$S(V)\ast G$, since it is finitely generated over $S(V)^{G}$. Thus
$\text{End}_{S(V)^{G}}\big(S(V)\big)$ is an order in its quotient
ring $\text{End}_{Q(S(V)^{G})}\big[Q(S(V))\big]$ equivalent to the
maximal order im$\big(j\circ gr(\eta)\big)$. Now the maximality of
im$\big(j\circ gr(\eta)\big)$ implies that im$\big(j\circ
gr(\eta)\big)=\text{End}_{S(V)^{G}}\big(S(V)\big)$.
\end{proof}

\begin{prop}\label{grHA:spherical-centre}
Assume that $G=S$. Then $Z(eA_{\textbf{t},\textbf{c}}e)\cong
Z(A_{\textbf{t},\textbf{c}})$ as $\C$-algebras.
\end{prop}

\begin{proof} We adapt the proof in \cite[Theorem
3.1]{EG02} very slightly and mention it for completeness.

Define a $\C$-algebra map
$\psi:Z(A_{\textbf{t},\textbf{c}})\rightarrow
Z(eA_{\textbf{t},\textbf{c}}e)$ by $z\mapsto ze=eze$ for $z\in
Z(A_{\textbf{t},\textbf{c}})$. We want to construct an inverse
algebra map to $\psi$ denoted by
$\varphi:Z(eA_{\textbf{t},\textbf{c}}e)\rightarrow
Z(A_{\textbf{t},\textbf{c}})$. Say $eae\in
Z(eA_{\textbf{t},\textbf{c}}e)$ and let $r_{eae}$ be right
multiplication by $eae$. Then $r_{eae}$ is an element of
$\text{End}_{eA_{\textbf{t},\textbf{c}}e}(A_{\textbf{t},\textbf{c}}e)$.
By the isomorphism
$\eta:A_{\textbf{t},\textbf{c}}\rightarrow\text{End}_{eA_{\textbf{t},\textbf{c}}e}(A_{\textbf{t},\textbf{c}}e)$
of the previous lemma we have 
$r_{eae}=\eta(x(a))=l_{x(a)}$ for some $x(a)\in
A_{\textbf{t},\textbf{c}}$, where $l_{x(a)}$ denotes left
multiplication by $x(a)$. Since left multiplication commutes with
right multiplication in
$\text{End}_{eA_{\textbf{t},\textbf{c}}e}(A_{\textbf{t},\textbf{c}}e)$,
$r_{eae}=l_{x(a)}$ is central in
$\text{End}_{eA_{\textbf{t},\textbf{c}}e}(A_{\textbf{t},\textbf{c}}e)$.
Now the isomorphism $\eta$ implies that $x(a)\in
Z(A_{\textbf{t},\textbf{c}})$. Thus define $\varphi: eae\mapsto
x(a)$. This is an algebra map because the isomorphism $\eta$ is an
algebra map. It remains to show that $\psi$ and $\varphi$ are
inverse to each other. We have $\varphi\circ\psi:z\mapsto eze\mapsto
x(z)$. As $z$ is central, we have $r_{eze}=r_{z}$ and $r_{z}=l_{z}$.
This implies that $l_{x(z)}=l_{z}$, that is $\eta(x(z))=\eta(z)$,
which implies $x(z)=z$, because $\eta$ is an isomorphism. On the
other hand $\psi\circ\varphi:eae\mapsto x(a)\mapsto ex(a)e$. For all
$y$ in $Z(eA_{\textbf{t},\textbf{c}}e)$ we have $l_{x(a)}(y)\cdot
e=r_{eae}(y)\cdot e$. But $l_{x(a)}(y)\cdot e=x(a)\cdot y\cdot
e=y\cdot x(a)\cdot e=y\cdot ex(a)e$, because $x(a)\in
Z(A_{\textbf{t},\textbf{c}})$, and $r_{eae}(y)\cdot e=y\cdot
eae\cdot e=y\cdot eae$. Thus $y\cdot ex(a)e-y\cdot eae =0$ and
$y[(ex(a)e-eae]=0$. Since $eA_{\textbf{t},\textbf{c}}e$ does not
contain zero divisors by Lemma \ref{grHA:spherical-domain}, this
implies $ex(a)e=eae$ as required.
\end{proof}
\section{Preliminary results}\label{grHA-preliminary}

Recall from Section \ref{grHA-definition} that
$A_{\textbf{t},\textbf{c}}(S)$ is defined as the subalgebra of
$A_{\textbf{t},\textbf{c}}$ constructed with the subgroup $S$ of
$G$. We can reduce to the case $G=S$ without loss of generality for
our main theorem because of the following result:

\begin{lem} \label{grHA:reduction} If $A_{\textbf{t},\textbf{c}}(S)$
is a finitely generated module over its centre
$Z(A_{\textbf{t},\textbf{c}}(S))$, then $A_{\textbf{t},\textbf{c}}$
is a finitely generated module over its centre and a PI algebra.
\end{lem}

\begin{proof} 
In Proposition \ref{grHA:firstproperties} we saw that
$A_{\textbf{t},\textbf{c}}\cong A_{\textbf{t},\textbf{c}}(S)\ast'
G/S$. This implies that $A_{\textbf{t},\textbf{c}}$ is finitely
generated over $A_{\textbf{t},\textbf{c}}(S)$ and that
$Z(A_{\textbf{t},\textbf{c}})\supseteq
\big[Z(A_{\textbf{t},\textbf{c}}(S))\big]^{G/S}$. 
It now suffices to show that $A_{\textbf{t},\textbf{c}}(S)$ is
finitely generated over
$\big[Z(A_{\textbf{t},\textbf{c}}(S))\big]^{G/S}$. But by the
initial assumption it only remains to show that
$Z(A_{\textbf{t},\textbf{c}}(S))$ is finitely generated over
$\big[Z(A_{\textbf{t},\textbf{c}}(S))\big]^{G/S}$. We have
$\C\subseteq Z(A_{\textbf{t},\textbf{c}}(S))\subseteq
A_{\textbf{t},\textbf{c}}(S)$, and $A_{\textbf{t},\textbf{c}}(S)$ is
an affine $\C$-algebra, which is a finite
$Z(A_{\textbf{t},\textbf{c}}(S))$-module. Thus the Artin-Tate lemma,
see \cite[Lemma 13.9.11]{MR87}, implies that
$Z(A_{\textbf{t},\textbf{c}}(S))$ is an affine $\C$-algebra as well.
But $G/S$ acts as a group of automorphisms on
$Z(A_{\textbf{t},\textbf{c}}(S))$ and we can use the Hilbert-Noether
theorem, see \cite[Theorem 1.3.1]{Ben93}, to deduce that
$Z(A_{\textbf{t},\textbf{c}}(S))$ is a finite
$\big[Z(A_{\textbf{t},\textbf{c}}(S))\big]^{G/S}$-module. Now
$A_{\textbf{t},\textbf{c}}$ is finitely generated over a commutative
subalgebra and hence a PI algebra by \cite[Corollary 13.1.13]{MR87}.
\end{proof}


Conversely, if $A_{\textbf{t},\textbf{c}}$ is a PI algebra then its
subalgebra $A_{\textbf{t},\textbf{c}}(S)$ is also a PI algebra, see
\cite[Lemma 13.1.7]{MR87}. In general, this does not imply that
$A_{\textbf{t},\textbf{c}}(S)$ is finitely generated over its centre
$Z(A_{\textbf{t},\textbf{c}}(S))$. However, we will derive this
implication as a consequence of our main theorem in the last section.\\

Throughout the remainder of this section we assume that $G=S$.\\

Recall that a \emph{Poisson bracket} on a commutative $\C$-algebra,
say $S(V)^{G}$, is a bilinear map $\{-,-\}:S(V)^{G}\times
S(V)^{G}\rightarrow S(V)^{G}$ such that $S(V)^{G}$ is a Lie algebra
under the bracket $\{-,-\}$ and the
Leibniz identity holds. 
In particular, $\{-,-\}$ satisfies the Jacobi identity. Moreover, a
Poisson bracket on $S(V)^{G}$ can be identified with an element of
$\text{Hom}_{S(V)^{G}}\big(\bigwedge^{2}D_{S(V)^{G}/\C},S(V)^{G}\big)$,
where $D_{S(V)^{G}/\C}$ denotes the module of K\"{a}hler
differentials of $S(V)^{G}$ over $\C$. 
The module $D_{S(V)^{G}/\C}$ is an $S(V)^{G}$-module and the
generators of $D_{S(V)^{G}/\C}$ are denoted by $dp$ for $p\in
S(V)^{G}$. The identification of a bracket $\{-,-\}$ with $\alpha\in
\text{Hom}_{S(V)^{G}}\big(\bigwedge^{2}D_{S(V)^{G}/\C},S(V)^{G}\big)$
is as follows: for $p,p'\in S(V)^{G}$, $\{p,p'\}\mapsto \big(\alpha:
dp\wedge dp'\mapsto \{p,p'\}\big)$. The Jacobi identity on $\{-,-\}$
imposes a relation on the map $\alpha$.

The algebra $S(V)^{G}$ is graded using the usual grading on $S(V)$.
Denote the $i$th graded part of $S(V)^{G}$ by $S^{i}(V)^{G}$ and
observe that $S^{i}(V)^{G}=0$ for $i<0$. A Poisson bracket $\{-,-\}$
on $S(V)^{G}$ is said to have degree $d$ if
$\{-,-\}:S^{i}(V)^{G}\times S^{j}(V)^{G}\rightarrow
S^{i+j+d}(V)^{G}$. Note that each element $\omega$ of
$((\bigwedge^{2}V)^{\ast})^{G}$ induces a Poisson bracket on $S(V)$
by extending $\omega$ linearly and using the Leibniz rule. Let us
denote this bracket by $\{-,-\}_{\omega}$. The fact that $\omega$ is
$G$-invariant forces $\{-,-\}_{\omega}$ to be a $G$-invariant
bracket on $S(V)$ as well. Thus $\{-,-\}_{\omega}$ restricts to a
Poisson bracket on $S(V)^{G}$. Furthermore, the bracket
$\{-,-\}_{\omega}$ on $S(V)^{G}$ has degree $-2$.

\begin{lem}\label{grHA:Poisson} Any Poisson bracket on $S(V)^{G}$ of
degree $-2$ is induced by an element of
$((\bigwedge^{2}V)^{\ast})^{G}$. Any Poisson bracket of degree less
than $-2$ is zero.
\end{lem}

\begin{proof} We proceed as in the proof of \cite[Lemma
2.23]{EG02}, but do not assume that the vector space $V$ comes
equipped with a symplectic form. Full details are given for the
convenience of the reader. 

Let $\{-,-\}$ denote a Poisson bracket on
$S(V)^{G}=\mathcal{O}(V^{\ast}/G)$ of degree $d$. In the proof of
this lemma we will extend the bracket $\{-,-\}$ on $S(V)^{G}$ to a
$G$-invariant Poisson bracket of degree $d$ on
$S(V)=\mathcal{O}(V^{\ast})$, denoted by $\{\{-,-\}\}$. We are then
able to prove that such a bracket on $S(V)$ is zero for $d<-2$ and
that it has to be induced by an element of
$((\bigwedge^{2}V)^{\ast})^{G}$ for $d=-2$.

In order to construct the bracket $\{\{-,-\}\}$, we first pick a
smooth open subset of $V^{\ast}/G$ as follows. Let $Y$ be the set of
points in $V^{\ast}$ that are fixed by some nontrivial element of
$G$, so $Y=\cup_{g\in G, g\neq 1} (V^{\ast})^{g}$. Note that
$(V^{\ast})^{g}$ is the zero set of the ideal $I_{g}\lhd S(V)$ given
by $I_{g}=\langle gv-v\,:\,v\in V\rangle$. Since the action of $G$
on $V$ is faithful, $I_{g}\neq 0$. Hence $(V^{\ast})^{g}$ is a
proper closed subset of $V^{\ast}$ for all $g\in G$. Then $Y$ is the
zero set of $I:=\cap_{g\in G, g\neq 1} I_{g}$ and a proper closed
subset of $V^{\ast}$. Therefore, the open set $X:= V^{\ast}\setminus
Y$ is a quasi-affine variety. Furthermore, the action of $G$ on $X$
is free in a set-theoretic sense, that is for all $x\in X$ the
stabiliser of $x$ in $G$, denoted by $G_{x}$, is trivial. Now
\cite[Proposition 4.12]{Dre04} says that the quotient map
$\pi:V^{\ast}\rightarrow V^{\ast}/G$ is \'etale at $x\in V^{\ast}$
if and only if $G_{x}$ is trivial. A consequence of $\pi$ being
\'etale at all $x\in X$ is an isomorphism between the completions of
the local rings $\mathcal{O}(V^{\ast}/G)_{\pi(x)}$ and
$\mathcal{O}(V^{\ast})_{x}$, that is
$\widehat{\mathcal{O}(V^{\ast}/G)}_{\pi(x)}\cong\widehat{\mathcal{O}(V^{\ast})}_{x}$,
for all $x\in X$, see \cite[Proposition 4.2]{Dre04}. Since
$V^{\ast}$ is a smooth variety, the local ring
$\mathcal{O}(V^{\ast})_{x}$ is regular for all $x\in V^{\ast}$,
which implies that $\widehat{\mathcal{O}(V^{\ast})}_{x}$ is regular
for all $x\in V^{\ast}$, see \cite[Theorem I.5.1, Theorem
I.5.4A]{Har77}. Thus, using the same results, we deduce that
$\mathcal{O}(V^{\ast}/G)_{\pi(x)}$ is regular for all $\pi(x)\in
V^{\ast}/G$ such that $x\in X$. But $\pi|_{X}:X\rightarrow X/G$ is
surjective, hence $X/G$ is a smooth variety.

Now take the given Poisson bracket $\{-,-\}$ on $S(V)^{G}$ of degree
$d$. For any open subset $U$ of $V^{\ast}/G$, $\{-,-\}$ defines a
map $\mathcal{O}(U)\times\mathcal{O}(U)\rightarrow\mathcal{O}(U)$.
Since the quotient map $\pi$ is closed, it takes the open subset
$X\subseteq V^{\ast}$ to an open subset $X/G$ of $V^{\ast}/G$. Hence
the bracket $\{-,-\}$ restricts to a 
Poisson bracket of degree $d$ on the sheaf of regular functions
$\mathcal{O}_{X/G}$ of the smooth variety $X/G$. 
We now observe that we can lift this bracket on $\mathcal{O}_{X/G}$
to a $G$-invariant Poisson bracket of degree $d$ on
$\mathcal{O}_{X}$. The reason for this is that the action of $G$ on
$X$ is free in a set-theoretic sense, so the quotient map
$\pi|_{X}:X\rightarrow X/G$ is not only \'etale but also a Galois
cover, see \cite[Definition 6.1]{Mil98}. We saw that a Poisson
bracket $\{-,-\}$ on $\mathcal{O}_{X/G}$ can be identified with an
element of
$\text{Hom}_{\mathcal{O}_{X/G}}\big(\bigwedge^{2}D_{\mathcal{O}_{X/G}/\C},\mathcal{O}_{X/G}\big)$,
which is the set of global sections of the second exterior power of
the tangent sheaf on $X/G$, see \cite[Definition, p.180]{Har77} for
the definition of a tangent sheaf. Now the theory on \'etale sheafs
and Galois coverings, as outlined in \cite[Section 6]{Mil98}, allows
one to identify
$\text{Hom}_{\mathcal{O}_{X/G}}\big(\bigwedge^{2}D_{\mathcal{O}_{X/G}/\C},\mathcal{O}_{X/G}\big)$
with
$\Big[\text{Hom}_{\mathcal{O}_{X}}\big(\bigwedge^{2}D_{\mathcal{O}_{X}/\C},\mathcal{O}_{X}\big)\Big]^{G}$.
Let us denote the resulting $G$-invariant Poisson bracket of degree
$d$ on $\mathcal{O}_{X}$ by $\{-,-\}_{X}$.

Recall that $X=V^{\ast}\setminus Y$. The next step is to extend the
bracket $\{-,-\}_{X}$ on $\mathcal{O}_{X}$ to a $G$-invariant
Poisson bracket on $\mathcal{O}(V^{\ast})=S(V)$. Since the group
$G=S\subseteq SL(V)$ does not contain reflections, each non-identity
element in $G$ has at least two eigenvalues different from $1$. This
implies that the codimension of $V^{g}$ is at least $2$ for all
$g\in G$, which translates into the corresponding ideal $I_{g}$
having height at least $2$. Hence the height of $I$ and the
codimension of $Y$ in $V$ is at least $2$ as well, see \cite[Section
II.1.3]{Kun85}. This enables us to apply \cite[Theorem
1.5.14]{FSR05} to extend a regular element
$x\in\mathcal{O}(X)=\mathcal{O}_{X}(X)$ to a regular element
$\tilde{x}$ in $\mathcal{O}(V^{\ast})$ such that $\tilde{x}|_{X}=x$.
Note that this is well-defined: say $x=x'\in\mathcal{O}(X)$. Then we
must have $\tilde{x}=\tilde{x'}\in\mathcal{O}(V^{\ast})$, because
$\tilde{x}$ and $\tilde{x'}$ agree on the non-empty and therefore
dense open subset $X$ of $V^{\ast}$. Furthermore, the map
$\mathcal{O}(V^{\ast})\twoheadrightarrow \mathcal{O}(X)$ given by
restriction is a surjection. Thus we construct a Poisson bracket on
$\mathcal{O}(V^{\ast})=S(V)$ denoted by $\{\{-,-\}\}$ as follows:
for $\tilde{x},\tilde{x'}\in \mathcal{O}(V^{\ast})$, define
$\{\{\tilde{x},\tilde{x'}\}\}:=\widetilde{\{x,x'\}}_{X}$. Since the
bracket $\{-,-\}_{X}$ is $G$-invariant and of degree $d$, the new
Poisson bracket $\{\{-,-\}\}$ is $G$-invariant and of degree $d$ as
well.

The bracket $\{\{-,-\}\}$ on $S(V)$ corresponds to a $G$-invariant
element of degree $d$ in
$\text{Hom}_{S(V)}\big(\bigwedge^{2}D_{S(V)/\C},S(V)\big)$. On the
other hand we have,
$\bigwedge^{2}\text{Hom}_{S(V)}\big(D_{S(V)/\C},S(V)\big) \cong
\text{Hom}_{S(V)}\big(\bigwedge^{2}D_{S(V)/\C},S(V)\big)$. 
Furthermore, $\text{Hom}_{S(V)}\big(D_{S(V)/\C},S(V)\big) \cong
\text{Der}_{\C}(S(V))$, the latter being the algebra of
$\C$-derivations on $S(V)$, see \cite[Proposition 15.1.10]{MR87}. In
summary,
$\text{Hom}_{S(V)}\big(\bigwedge^{2}D_{S(V)/\C},S(V)\big)\cong\bigwedge^{2}\text{Der}_{\C}(S(V))$.
It is now easy to see that an element of degree $d<-2$ in
$\bigwedge^{2}\text{Der}_{\C}(S(V))$ is zero. Hence if the degree of
$\{\{-,-\}\}$ is less than $-2$, then the bracket $\{\{-,-\}\}$ is
zero. Furthermore, an element of degree $-2$ in
$\bigwedge^{2}\text{Der}_{\C}(S(V))$ must be an element of
$\bigwedge^{2} V^{\ast}$. If we assume in addition that this element
is $G$-invariant, then it must be an element of $(\bigwedge^{2}
V^{\ast})^{G}\cong ((\bigwedge^{2}V)^{\ast})^{G}$.
\end{proof}

\begin{lem}\label{grHA:regrep} Suppose the parameters in $\textbf{t}$ and
$\textbf{c}$ are such that $eA_{\textbf{t},\textbf{c}}e$ is
commutative. Let \emph{MaxSpec}$(eA_{\textbf{t},\textbf{c}}e)$
denote the set of maximal ideals of $eA_{\textbf{t},\textbf{c}}e$.
Then there exists a non-empty Zariski-open subset $\mathcal{M}$ of
\emph{MaxSpec}$(eA_{\textbf{t},\textbf{c}}e)$ such that, if
$\mathfrak{m}\in\mathcal{M}$ and if we let
$T_{\mathfrak{m}}:=A_{\textbf{t},\textbf{c}}e\otimes_{eA_{\textbf{t},\textbf{c}}e}(eA_{\textbf{t},\textbf{c}}e/\mathfrak{m})$
denote the corresponding induced $A_{\textbf{t},\textbf{c}}$-module,
then $T_{\mathfrak{m}}\cong \C G$ as $G$-module.
\end{lem}

\begin{proof} The proof is the same as the one for \cite[Lemma
2.24]{EG02}.
\end{proof}
\section{Proof of the main theorem}\label{grHA-main}

Let $A_{\textbf{t},\textbf{c}}$ be a graded Hecke algebra as defined
in Section \ref{grHA-definition} and recall that
$A_{\textbf{t},\textbf{c}}$ is completely determined by the values
chosen for the parameters $\{t_{i}\,|\,i=1,\ldots, N\}$ and
$\{c_{s}\,|\,s\in \mathcal{S}'\}$, see Theorem \ref{grHA:def}. We
continue to assume for now that $G$ is generated by the elements in
$\mathcal{S}'$, hence $G=S$.

\begin{rem} It is probably possible to obtain the following result
for all finite groups $G\subseteq GL(V)$, that is to drop the
assumption $G=S$. However, it is not trivial to prove this and we do
not need this version for our purposes.
\end{rem}

\begin{thrm}\label{grHA:spherical-comm} Assume $G=S$. Then $eA_{\textbf{t},\textbf{c}}e$ is commutative if and
only if $t_{i}=0$ for all $i=1,\ldots,N$.
\end{thrm}

\begin{proof} The proof of this theorem uses the deformation theory that we introduced in Section
\ref{grHA-deformation} and a strategy similar to the one in the
proof of \cite[Theorem 1.6]{EG02}.

Since $A_{\textbf{t},\textbf{c}}$ is a PBW deformation of
$R=S(V)\ast G$, as seen in Corollary \ref{grHA:deform}, there exists
a graded deformation $(R_{h}, \ast)$ of $R$ such that
$R_{h}/(h-1)R_{h}\cong A_{\textbf{t},\textbf{c}}$. In order to
describe such a deformation $R_{h}$ explicitly we introduce the
auxiliary variable $h$ and set $T(V)[h]:=T(V)\otimes \C[h]$. Let the
degree of $h$ be $1$ and assume that the group $G$ acts trivially on
$h$. Define
\[R_{h}:=\big(T(V)[h]\ast G\big)/\langle [v,w]-\kappa(v,w)h^{2}\,:\,v,w\in V\rangle.\]

The algebra $R_{h}$ is indeed a graded deformation of $R$. 
Namely, since the relation $[v,w]=\kappa(v,w)h^{2}$ is now
homogeneous,
$R_{h}$ is an associative unital graded algebra. 
It is easy to see that $R_{h}/hR_{h}=S(V)\ast G=R$ and that
$R_{h}/(h-1)R_{h}=A_{\textbf{t},\textbf{c}}$. If we pick a vector
space basis $v_{1},\ldots, v_{n}$ of $V$, we obtain a vector space
basis of $S(V)$ consisting of ordered monomials in the $v_{i}$. Now
we can think of $p\in S(V)$ as an element of $T(V)$. We can then use
the projection $T(V)\otimes \C G\rightarrow R_{h}$ to obtain an
epimorphism of $\C$-vector spaces $\pi:\big(S(V)\otimes\C
G\big)[h]\twoheadrightarrow R_{h}$ given by
$\sum_{i=0}^{m}p_{i}h^{i}\mapsto \sum_{i=0}^{m}p_{i}h^{i}$, where
$p_{i}\in S(V)\otimes\C G$, that is $p_{i}=\sum_{g\in G}p_{i,g}g$
and each $p_{i,g}$ is a linear combination of ordered monomials in
the $v_{i}$. We want to show that $\pi$ is an isomorphism, hence
that the underlying vector space of $R_{h}$ is $R\otimes\C[h]$. Thus
we need to prove that $\pi\big(\sum_{i=0}^{m}p_{i}h^{i}\big)=0$
implies $p_{i}=0$ for all $i=0,\ldots,m$. The map $\pi$ is a
homogeneous map of degree zero. Hence we can assume without loss of
generality that $\sum_{i=0}^{m}p_{i}h^{i}$ is a homogeneous element
of degree $k$. So $p_{i}\in S(V)\otimes\C G$ has degree $k-i$.
Denote the projection $R_{h}\rightarrow
R_{h}/(h-1)R_{h}=A_{\textbf{t},\textbf{c}}$ by $\varrho$. If
$\pi\big(\sum_{i=0}^{m}p_{i}h^{i}\big)=0$, then
$\varrho\big(\pi\big(\sum_{i=0}^{m}p_{i}h^{i}\big)\big)=\sum_{i=0}^{m}p_{i}+(h-1)R_{h}=0$.
This holds if and only if $\sum_{i=0}^{m}p_{i}\in(h-1)R_{h}$ which
is the case if and only if $\sum_{i=0}^{m}p_{i}=0\in R_{h}$. But the
elements $p_{i}\in R_{h}$ have distinct degrees, which means that we
must have $p_{i}=0$ for all $i=0,\ldots,m$ as required.

The multiplication $\ast$ in $R_{h}$ is given by the multiplication
in $T(V)\ast G$ and the additional relations $v\ast w-w\ast
v=\kappa(v,w)h^{2}$ for all $v,w\in V$, which are extended by
$\C[h]$-linearity. We have seen during the proof of Corollary
\ref{grHA:deform} that $\kappa(-,-)$ is a $\C G$-bimodule map, which
makes $\ast$ into a $\C G[h]$-bimodule map. Thus $R_{h}$ is indeed a
graded $\C G[h]$-bimodule.
We can express multiplication in $R_{h}$ in terms of $\C G$-bimodule
maps $\mu_{i}: R\times R\rightarrow R$ of degree $-i$. At the end of
Section \ref{grHA-deformation} 
we observed that $\kappa(v,w)=\mu_{2}(v,w)-\mu_{2}(w,v)$ and
$\mu_{1}(v,w)-\mu_{1}(w,v)=0$, for all $v,w\in V$. Furthermore,
$p\ast p'-p'\ast p=\mu_{2}(p,p')h^2-\mu_{2}(p',p)h^{2}$ for all
$p,p'\in S(V)$.

Let us form the 
spherical subalgebra $eR_{h}e$ of $R_{h}$. Clearly, 
$(eR_{h}e,\ast)$ is a graded deformation of $e(S(V)\ast G)e\cong
S(V)^{G}$, that is $eR_{h}e/heR_{h}e=e(S(V)\ast G)e$. This is
because we chose the maps $\mu_{i}$ to be $\C G$-invariant. As
vector spaces, $eR_{h}e\cong S(V)^{G}[h]$. Also,
$eR_{h}e/(h-1)eR_{h}e=eA_{\textbf{t},\textbf{c}}e$. Given $p,p'\in
S(V)^{G}\cong eR_{h}e/heR_{h}e$ let $\tilde{p},\tilde{p'}$ denote
lifts of these elements to $eR_{h}e$. 
We define a Poisson bracket $\{-,-\}$ on $S(V)^{G}$ by
$\{p,p'\}:=h^{-2}(\tilde{p}\ast\tilde{p'}-\tilde{p'}\ast\tilde{p})\,\text{mod}(heR_{h}e)$.
It is easy to check that this indeed defines a Poisson bracket and
that $\{p,p'\}=\mu_{2}(p,p')-\mu_{2}(p',p)$, for all $p,p'\in
S(V)^{G}$.
We claim that
\begin{equation}\label{grHA:spherical-comm-poisson}
eA_{\textbf{t},\textbf{c}}e\,\text{commutative}\,\Leftrightarrow\,
eR_{h}e\,\text{commutative}\,\Leftrightarrow\,\{-,-\}\equiv 0
\end{equation}
Let us first show the equivalence on the right hand side. From the
last description of the Poisson bracket it becomes obvious that, if
$eR_{h}e$ is commutative, then $\{-,-\}\equiv 0$. Conversely, if
$\{p,p'\}=0$ for all $p,p'\in S(V)^{G}$, then
$\mu_{2}(p,p')=\mu_{2}(p',p)$ for all $p,p'\in S(V)^{G}$. Since
$eR_{h}e$ is a deformation of $S(V)^{G}$, the multiplication $\ast$
in $eR_{h}e$ is determined by the multiplication $S(V)^{G}\ast
S(V)^{G}\subset eR_{h}e\ast eR_{h}e$ and then extended by
$\C[h]$-linearity. But we now have, for all $p,p'\in S(V)^{G}$,
$p\ast p'-p'\ast p=[\mu_{2}(p,p')-\mu_{2}(p',p)]h^{2}=0$. Hence
$eR_{h}e$ is commutative. For the equivalence on the left hand side
we observe that, if $eR_{h}e$ is commutative, the factor algebra
$eR_{h}e/(h-1)eR_{h}e=eA_{\textbf{t},\textbf{c}}e$ is certainly also
commutative. Conversely, assume $eR_{h}e/(h-1)eR_{h}e$ is
commutative, but $eR_{h}e$ is not. Then, by the above, the Poisson
bracket is nonzero and so there exist $p,p'\in eR_{h}e/heR_{h}e$
such that $\{p,p'\}=f\neq 0$. Choose representatives
$\tilde{p},\tilde{p'}\in eR_{h}e$ of $p,p'$. We can assume without
loss of generality that $\tilde{p},\tilde{p'}$ are homogeneous
elements of $eR_{h}e$. Then $h^{-2}(\tilde{p}\ast
\tilde{p'}-\tilde{p'}\ast\tilde{p})=\tilde{f}$ such that
$f\equiv\tilde{f}\,\text{mod}(heR_{h}e)$ and $\tilde{f}$ is a
nonzero homogeneous element of $eR_{h}e$. Now consider
$[\tilde{p},\tilde{p'}]$ mod $\big((h-1)eR_{h}e\big)$. Since
$eR_{h}e/(h-1)eR_{h}e=eA_{\textbf{t},\textbf{c}}e$ is assumed to be
commutative, $[\tilde{p},\tilde{p'}]\equiv
0\,\text{mod}\big((h-1)eR_{h}e\big)$. But
$[\tilde{p},\tilde{p'}]=h^{2}\tilde{f}$, i.e
$[\tilde{p},\tilde{p'}]\equiv
\tilde{f}\,\text{mod}\big((h-1)eR_{h}e\big)$. Thus
$\tilde{f}\in\big((h-1)eR_{h}e\big)$, which means that $\tilde{f}$
is divisible by $(h-1)$. We conclude that $\tilde{f}$ is not
homogeneous, a contradiction. 

It now remains to prove that the Poisson bracket $\{-,-\}$ on
$S(V)^{G}$ vanishes if and only if $t_{i}=0$ for all $i=1,\ldots,
N$. 
Since the degree of the map $\mu_{2}$ is $-2$, the degree of the
Poisson bracket is also $-2$. Hence Lemma \ref{grHA:Poisson} implies
that the bracket is induced by some element $\omega$ of
$((\bigwedge^{2}V)^{\ast})^{G}$. In terms of the basis
$\{b_{1},\ldots, b_{N}\}$ of $((\bigwedge^{2}V)^{\ast})^{G}$ write
$\omega=\sum_{i=1}^{N}\lambda_{i}b_{i}$, for some $\lambda_{i}\in
\C$. Let $\{-,-\}_{i}$ denote the Poisson bracket induced by
$b_{i}$. From the explanations preceeding Lemma \ref{grHA:Poisson}
it is easy to see that we must have
$\{-,-\}=\sum_{i=1}^{N}\lambda_{i}\{-,-\}_{i}$. Furthermore, at the
end of Section \ref{grHA-deformation} we observed that the
difference $\mu_{2}(p,p')-\mu_{2}(p',p)$ depends linearly on the
parameters $\textbf{t}$ and $\textbf{c}$ for all $p,p'\in S(V)^{G}$.
Thus the Poisson bracket $\{-,-\}$ depends linearly on the
parameters $\textbf{t}$ and $\textbf{c}$, 
and so do the scalars $\lambda_{i}$. 
Let $f_{i}:\C^{N}\times\C^{|\mathcal{S}|}\rightarrow\C$ denote
linear functions and write
$\{-,-\}=\sum_{i=1}^{N}f_{i}(\textbf{t},\textbf{c})\{-,-\}_{i}$. Now
the Poisson bracket vanishes if and only if
$f_{i}(\textbf{t},\textbf{c})=0$ for all $i=1,\ldots,N$, since the
brackets $\{-,-\}_{i}$ are linearly independent by the linear
independency of the basis elements
$b_{i}\in((\bigwedge^{2}V)^{\ast})^{G}$. We need to show that this
is the case if and only if $t_{i}=0$ for all $i=1,\ldots,N$.

The equations $f_{i}(\textbf{t},\textbf{c})=0$, $i=1,\ldots,N$, form
a system of homogeneous linear equations of rank $r\leq N$. Thus the
solution space
$\mathcal{V}(f_{i})\subseteq\C^{N}\oplus\C^{|\mathcal{S}|}$ of these
equations has dimension $(N+|\mathcal{S}|)-r\geq
(N+|\mathcal{S}|)-N=|\mathcal{S}|$. On the other hand, the system of
linear equations given by $t_{i}=0$, $i=1,\ldots,N$, has rank $N$
and, therefore, its solution space, $\mathcal{V}(t_{i})$, is
$|\mathcal{S}|$-dimensional. Thus
dim$\big(\mathcal{V}(f_{i})\big)\geq\text{dim}\big(\mathcal{V}(t_{i})\big)$.
We will show that $\mathcal{V}(f_{i})\subseteq\mathcal{V}(t_{i})$
which implies the result by containment and equality of dimensions.
To show that $\mathcal{V}(f_{i})\subseteq\mathcal{V}(t_{i})$ we
assume that the parameters $\textbf{t},\textbf{c}$ are such that
$f_{i}(\textbf{t},\textbf{c})=0$ for all $i=1,\ldots,N$. Then the
Poisson bracket on $S(V)^{G}$ vanishes and
$eA_{\textbf{t},\textbf{c}}e$ is commutative. We can now use Lemma
\ref{grHA:regrep}. Let
$\mathfrak{m}\in\text{MaxSpec}(eA_{\textbf{t},\textbf{c}}e)$ be such
that the corresponding induced $A_{\textbf{t},\textbf{c}}$-module
$T_{\mathfrak{m}}=A_{\textbf{t},\textbf{c}}e\otimes_{eA_{\textbf{t},\textbf{c}}e}(eA_{\textbf{t},\textbf{c}}e/\mathfrak{m})$
is isomorphic to $\C G$ as $G$-module. In
$A_{\textbf{t},\textbf{c}}$ we have the relation $v\otimes
w-w\otimes v=\kappa(v,w)\in\C G$, for all $v,w\in V$. Now 
take traces on both sides of this equation: tr$(v\otimes w-w\otimes
v)=0$ and
tr$\big(\kappa(v,w)\big)=\text{tr}\Big(\big[\sum_{i=1}^{N}t_{i}b_{i}(v,w)\big]\cdot
id+\sum_{s\in \mathcal{S}'}c_{s}\Omega_{s}(v,w)
s\Big)=\big[\sum_{i=1}^{N}t_{i}b_{i}(v,w)\big]
\text{tr}(id)+\sum_{s\in
\mathcal{S}'}c_{s}\Omega_{s}(v,w)\text{tr}(s)$. But because
$T_{\mathfrak{m}}$ is isomorphic to the regular representation of
$G$ as a $G$-module, tr$(id)=|G|$ and tr$(s)=0$ for all $s\neq 1$.
This implies that
$\sum_{i=1}^{N}t_{i}b_{i}=0\in((\bigwedge^{2}V)^{\ast})^{G}$ and by
the linear independence of the $b_{i}$ we conclude that $t_{i}=0$
for all $i=1,\ldots,N$.
\end{proof}


\begin{cor}\label{grHA:fingen}
\ \\
\indent (i) $eA_{\textbf{0},\textbf{c}}e\cong
Z(A_{\textbf{0},\textbf{c}})$ as $\C$-algebras.

(ii) $gr Z(A_{\textbf{0},\textbf{c}})\cong S(V)^{G}$.

(iii) $A_{\textbf{0},\textbf{c}}$ is a finitely generated module
over $Z(A_{\textbf{0},\textbf{c}})$ and $A_{\textbf{0},\textbf{c}}$
is a PI-algebra.
\end{cor}

\begin{proof} (i) Follows from Proposition
\ref{grHA:spherical-centre} and the previous theorem.

(ii) In Proposition \ref{grHA:spherical-centre} we found an
isomorphism $\psi:Z(A_{\textbf{0},\textbf{c}})\rightarrow
Z(eA_{\textbf{0},\textbf{c}}e)=eA_{\textbf{0},\textbf{c}}e$ given by
$z\mapsto ze$, for $z\in Z(A_{\textbf{0},\textbf{c}})$. The map
$\psi$ is filtration preserving since $e\in F^{0}$. Thus we have
$\psi\big(F^{i}_{Z(A_{\textbf{0},\textbf{c}})}\big)\subseteq\psi(Z(A_{\textbf{0},\textbf{c}}))\cap
F^{i}_{eA_{\textbf{0},\textbf{c}}e}$ for all $i\geq 0$. But if $z\in
Z(A_{\textbf{0},\textbf{c}})$ and $\psi(z)=ze\in
F^{i}_{eA_{\textbf{0},\textbf{c}}e}$, then we can easily see that
$z\in F^{i}_{A_{\textbf{0},\textbf{c}}}\cap
Z(A_{\textbf{0},\textbf{c}})$, because
$F^{i}_{eA_{\textbf{0},\textbf{c}}e}=eF^{i}_{A_{\textbf{0},\textbf{c}}}e$.
Now the surjectivity of $\psi$ implies that
$\psi\big(F^{i}_{Z(A_{\textbf{0},\textbf{c}})}\big)=
F^{i}_{eA_{\textbf{0},\textbf{c}}e}$ for all $i\geq 0$. Hence $\psi$
is a strict map, see \cite[7.6.12]{MR87}. Then \cite[Corollary
7.6.14]{MR87} implies that the induced map $gr
Z(A_{\textbf{0},\textbf{c}})\mapsto gr(eA_{\textbf{0},\textbf{c}}e)$
is bijective. But $gr(eA_{\textbf{0},\textbf{c}}e)\cong S(V)^{G}$ as
we saw at the beginning of Section \ref{grHA-spherical}.

(iii) It is enough to show that $gr A_{\textbf{0},\textbf{c}}$ is
finitely generated over $gr Z(A_{\textbf{0},\textbf{c}})$, because
we can then use associated graded arguments as in Lemma
\ref{grHA:spherical-domain}. Denote the isomorphism $\gamma:S(V)\ast
G\rightarrow gr A_{\textbf{0},\textbf{c}}$. Since $S(V)\ast G$ is
finitely generated over $S(V)^{G}=Z(S(V)\ast G)$, $\gamma(S(V)\ast
G)= gr A_{\textbf{0},\textbf{c}}$ is finitely generated over
$\gamma\big(Z(S(V)\ast G)\big)=Z\big(\gamma(S(V)\ast G)\big)= Z(gr
A_{\textbf{0},\textbf{c}})$. Thus it remains to prove that $gr
Z(A_{\textbf{0},\textbf{c}})=Z(gr A_{\textbf{0},\textbf{c}})$.
We have the following maps: $gr
Z(A_{\textbf{0},\textbf{c}})\rightarrow gr(e
A_{\textbf{0},\textbf{c}}e)=e (gr A_{\textbf{0},\textbf{c}})e$ given
by $z\mapsto ze$ for all $z\in gr Z(A_{\textbf{0},\textbf{c}})$ as
seen in Part (ii) of this corollary. And a map $S(V)^{G}\rightarrow
e (gr A_{\textbf{0},\textbf{c}})e$ given by $p\mapsto \gamma(p)e$.
Both of these maps are isomorphisms as observed in Part (ii) of this
corollary and at the beginning of Section \ref{grHA-spherical}. Thus
for each $\gamma(p)e$ there exists a unique $z\in gr
Z(A_{\textbf{0},\textbf{c}})$ such that $ze=\gamma(p)e$. Since
$\gamma(S(V)^{G})=Z(gr A_{\textbf{0},\textbf{c}})$, we can now
define a map $gr Z(A_{\textbf{0},\textbf{c}})\rightarrow Z(gr
A_{\textbf{0},\textbf{c}})$ by $z\mapsto \gamma(p)$. It is obvious
that this map is bijective. Now $A_{\textbf{0},\textbf{c}}$ is
finitely generated over a commutative subalgebra and hence a PI
algebra by \cite[Corollary 13.1.13]{MR87}.
\end{proof}
%

\begin{cor} \label{grHA:main} Let $A_{\textbf{t},\textbf{c}}$ be a graded Hecke algebra. Assume $G=S$.
Then $A_{\textbf{t},\textbf{c}}$ is a PI algebra if and only if
$A_{\textbf{t},\textbf{c}}$ is a finitely generated module over its
centre if and only if $t_{i}=0$ for all $i\in\{1,\ldots,N\}$.
\end{cor}

\begin{proof} From Theorem \ref{grHA:spherical-comm} and the
subsequent corollary we know
\[ A_{\textbf{t},\textbf{c}}\,\text{is a PI algebra}\Leftarrow
A_{\textbf{t},\textbf{c}}\,\text{is a finite}\,
Z(A_{\textbf{t},\textbf{c}})\text{-module}\Leftarrow
t_{i}=0\,\forall i\in\{1,\ldots,N\}. \] Thus it remains to prove
that if $A_{\textbf{t},\textbf{c}}$ is a PI algebra then $t_{i}=0$
for all $i\in\{1,\ldots,N\}$. To reach a contradiction assume that
$t_{i}\neq 0$ for some $i=1,\ldots,N$. This implies that the form
$\Omega=\sum_{i=1}^{N}t_{i}b_{i}$ is a nonzero $G$-invariant
skew-symmetric form on $V$. We claim that in this situation there
exists a subalgebra of $A_{\textbf{t},\textbf{c}}$ which is a
symplectic reflection algebra. Existing results on symplectic
reflection algebras will provide us with the necessary
contradiction.

Let $U:=\{u\in V\,|\,\Omega(u,v)=0\:\text{for all}\,v\in V\}$, the
radical of $\Omega$. Then $U$ is a $G$-invariant subspace of $V$,
because the form $\Omega$ is $G$-invariant. By Maschke's theorem, we
can find a $G$-invariant complement $W$ such that 
$V=U\oplus W$. Take $v, v'\in V$. We can write $v=u+w$ and
$v'=u'+w'$ for some $u,u'\in U, w,w'\in W$. We have
$\Omega(v,v')=\Omega(u+w,u'+w')=\Omega(w,w')$. Therefore, the form
$\Omega$ is determined by its restriction to $W$, denoted by
$\Omega|_{W}:W\times W\rightarrow\C$. Moreover, by construction, the
form $\Omega|_{W}$ is not only a nonzero $G$-invariant
skew-symmetric form on $W$, but also non-degenerate. In other words,
$W$ is a symplectic vector space.

Let $G'$ denote the subgroup of $G$ which is generated by those
elements that are bireflections in their action on the subspace $W$.
It is clear that $G'$ is closed under conjugation by elements of
$G$. We claim that the elements in $W$ and the elements in $G'$
generate a subalgebra of $A_{\textbf{t},\textbf{c}}$ which is a
symplectic reflection algebra. Obviously $T(W)\ast G'$ is a
subalgebra of $T(V)\ast G$. In order to prove our claim we need to
examine the relations
\[\kappa(w,w')=\Omega(w,w')\,id+\sum_{s\in
\mathcal{S}'}c_{s}\Omega_{s}(w,w')\, s,\] for all $w,w'\in W$. In
particular, we need to show that $\Omega_{s}|_{W\times W}=0$ for all
elements $s\in \mathcal{S}'$ that are not bireflections in their
action on $W$. Indeed, take $s\in\mathcal{S}'$. Since
dim$\big(\text{im}(id-s)\big)=2$, we have
dim$\big(\text{im}(id-s)\cap W\big)\leq 2$. 
Assume that dim$\big(\text{im}(id-s)\cap W\big)=0$, then $s$ fixes
$W$. But by construction, see Section \ref{grHA-definition}, the
subspace $V^{s}=\text{ker}(id-s)$ lies in the radical of
$\Omega_{s}$. Thus we deduce $\Omega_{s}|_{W\times W}=0$ for this
situation. Assume that dim$\big(\text{im}(id-s)\cap W\big)=1$. Say
im$(id-s)\cap W=\C x$. Since $\Omega_{s}|_{W\times W}$ is a
skew-symmetric form, $\Omega_{s}|_{W\times W}(\lambda x,\mu x)=0$
for all $\lambda,\mu\in\C$. But $W=(\text{im}(id-s)\cap W)\oplus
(\text{ker}(id-s)\cap W)$ and ker$(id-s)$ is again in the radical of
$\Omega_{s}$. Thus in this situation we also have
$\Omega_{s}|_{W\times W}=0$.

Denote the subalgebra of $A_{\textbf{t},\textbf{c}}$ generated by
$W$ and $G'$ by $A(W,G')$. Note that the action of $G'$ on $W$ is
faithful. Namely, the decomposition $V=U\oplus W$ is $G'$-invariant.
Take a generator $s$ of $G'\subseteq G$. Then $s$ is a bireflection
on $V$, because $s\in G$, but $s$ is also a bireflection on $W$. We
deduce that dim$\big(\text{im}(id-s)\cap U\big)=0$. So the group
$G'$ acts trivially on $U$. Now if $g\in G'$ is such that
$g|_{W}=id$, then $g|_{V}=id$. But because $G\subseteq GL(V)$ acts
faithfully on $V$, this implies that $g=id$. Therefore
$G'\hookrightarrow GL(W)$ and the subalgebra $A(W,G')$ is a
symplectic reflection algebra.

Since $A(W,G')$ is a subalgebra of the PI algebra
$A_{\textbf{t},\textbf{c}}$, it is also a PI algebra, see
\cite[Lemma 13.1.7]{MR87}. In \cite[Proposition 7.2]{BG03} it is
shown that if $\Omega|_{W}\neq 0$, then the centre
of the symplectic reflection algebra $A(W,G')$ is just $\C$. 
The fact that $A(W,G')$ is also prime, see Proposition
\ref{grHA:noeth,prime,gldim}, together with \cite[Proposition
13.6.11]{MR87} now implies that $A(W,G')$ is a finite dimensional
$\C$-vector space. But this is a contradiction to the fact that
$A(W,G')\cong S(W)\otimes\C G'$ as a $\C$-vector space.
\end{proof}

We now drop the assumption $G=S$ and finish with the general result:

\begin{cor} Let $A_{\textbf{t},\textbf{c}}$ be a graded Hecke
algebra. Then $A_{\textbf{t},\textbf{c}}$ is a PI algebra if and
only if $A_{\textbf{t},\textbf{c}}$ is a finitely generated module
over its centre if and only if $t_{i}=0$ for all
$i\in\{1,\ldots,N\}$.
\end{cor}

\begin{proof} We have $S\lhd G$ and we denote the graded Hecke algebra
constructed with $S$ instead of $G$ by
$A_{\textbf{t},\textbf{c}}(S)$. We have the implications:
\[\xymatrix
{A_{\textbf{t},\textbf{c}}\hspace{0.7mm}\text{PI}\ar@{=>}[d] & \ & A_{\textbf{t},\textbf{c}}\hspace{0.7mm}\text{PI}\\
A_{\textbf{t},\textbf{c}}(S)\hspace{0.7mm}\text{PI}\ar@{=>}[r] &
A_{\textbf{t},\textbf{c}}(S)\hspace{0.7mm} \text{a finite}\,
Z(A_{\textbf{t},\textbf{c}}(S))-\text{module}\ar@{=>}[r] &
A_{\textbf{t},\textbf{c}}\hspace{0.7mm} \text{a finite}\,
Z(A_{\textbf{t},\textbf{c}})-\text{module}\ar@{=>}[u]}\] where the
vertical implications are \cite[Lemma 13.1.7]{MR87} and
\cite[Corollary 13.1.13]{MR87}. The horizontal implications are the
corollary above and Lemma \ref{grHA:reduction}. Thus we know now
that $A_{\textbf{t},\textbf{c}}(S)$ is a finite
$Z(A_{\textbf{t},\textbf{c}}(S))$-module $\Leftrightarrow$
$A_{\textbf{t},\textbf{c}}$ is a finite
$Z(A_{\textbf{t},\textbf{c}})$-\text{module}. But, by the corollary
above, $A_{\textbf{t},\textbf{c}}(S)$ is a finite
$Z(A_{\textbf{t},\textbf{c}}(S))$-module $\Leftrightarrow$ $t_{i}=0$
for all $i\in\{1,\ldots,N\}$.
\end{proof}

\bibliography{PaperBib}
\bibliographystyle{alpha}

\end{document}